\documentclass[final]{amsart}
\usepackage{amsmath,amsthm,amsfonts,amssymb,indentfirst,xspace}
\usepackage{mathtools}
\usepackage[notref,notcite]{showkeys}
\usepackage[ps2pdf]{hyperref}
\usepackage[initials,nobysame]{amsrefs}

\newcommand{\snseqns}{\eqref{e:X-def}--\eqref{e:X-idata}\xspace}
\newcommand{\navierstokes}{\eqref{e:ns1}--\eqref{e:ns2}\xspace}

\newcommand{\del}{\partial}
\newcommand{\lap}{\triangle}
\newcommand{\inv}{^{-1}}
\newcommand{\transpose}{^*}
\newcommand{\leqs}{\leqslant}
\newcommand{\geqs}{\geqslant}
\newcommand{\grad}{\nabla}
\newcommand{\gradt}{\grad\transpose}
\newcommand{\divergence}{\grad \cdot}
\newcommand{\curl}{\grad \times}

\newcommand{\I}{\mathcal{I}}
\newcommand{\Ifn}{I}
\newcommand{\Imatrix}{\mathbb{I}}

\renewcommand{\epsilon}{\varepsilon}
\renewcommand{\leq}{\leqslant}
\renewcommand{\geq}{\geqslant}

\newcommand{\E}{\boldsymbol{\mathrm{E}}}
\newcommand{\lhp}{\boldsymbol{\mathrm{P}}}
\newcommand{\as}{\text{a.s.}}

\newcommand{\R}{\mathbb{R}}
\newcommand{\holderspace}[2]{\ensuremath{C^{\ifx0#1{#2}\else{#1,#2}\fi}}}
\newcommand{\cspace}[1]{\ensuremath{C^{#1}}}

\newif\iftextstyle
\textstyletrue
\everydisplay\expandafter{\the\everydisplay\textstylefalse}

\newcommand{\abs}[1]{{\iftextstyle\lvert#1\rvert\else\left\lvert#1\right\rvert\fi}}
\newcommand{\norm}[1]{{\iftextstyle\lVert#1\rVert\else\left\lVert#1\right\rVert\fi}}

\newcommand{\hnorm}[3]{\norm{#1}_{\ifx0#2{#3}\else{#2,#3}\fi}}
\newcommand{\hseminorm}[2]{\abs{#1}_{#2}}

\newcommand{\cnorm}[2]{\norm{#1}_{C^{#2}}}

\newcommand{\lpnorm}[2]{\norm{#1}_{L^{#2}}}

\newcommand{\at}[1]{{\bigr\rvert_{#1}}}

\numberwithin{equation}{section}
\allowdisplaybreaks

\newtheorem{theorem}{Theorem}[section]
\newtheorem*{theorem*}{Theorem}
\newtheorem{lemma}[theorem]{Lemma}
\newtheorem{proposition}[theorem]{Proposition}

\theoremstyle{definition}

\theoremstyle{remark}
\newtheorem{remark}[theorem]{Remark}
\newtheorem*{remark*}{Remark}

\begin{document}
\title[Global existence for Navier-Stokes with small Reynolds number.]{A stochastic Lagrangian proof of global existence of the Navier-Stokes equations for flows with small Reynolds number.}
\author{Gautam Iyer}
\address{%
Department of Mathematics\\
The University of Chicago\\
Chicago, Illinois 60637}
\email{gautam@math.uchicago.edu}
\begin{abstract}
We consider the incompressible Navier-Stokes equations with spatially periodic boundary conditions. If the Reynolds number is small enough we provide an elementary short proof of the existence of global in time H\"older continuous solutions. Our proof is based on the stochastic Lagrangian formulation of the Navier-Stokes equations, and works in both the two and three dimensional situation.
\end{abstract}
\maketitle
\section{Introduction}
The Navier-Stokes equations
\begin{gather}
\label{e:ns1} \del_t u + (u \cdot \grad) u - \nu \lap u + \grad p = 0\\
\label{e:ns2} \divergence u = 0
\end{gather}
describe the evolution of the velocity field of an incompressible fluid with kinematic viscosity $\nu > 0$. One of the (still open) million dollar problems posed by the Clay Institute \cite{clayproblem} is to show that given a smooth initial data $u_0$ the solution to \navierstokes in three dimensions remains smooth for all time.

In two dimensions, the long time existence of \navierstokes is well known \cite{const-book}. In three or higher dimensions, long time existence is known provided a smallness condition is imposed on the initial data (see for example \cites{bmo-1} for a criterion which in some sense is the most general smallness condition). Recently Chemin and Gallagher found a (non-linear) criterion on the initial data that guarantees global existence of \navierstokes, and does not reduce to a smallness criterion in $\text{BMO}\inv$.

In this paper we prove global existence of \navierstokes provided our initial data has small H\"older norm. Though global existence under our assumptions can be deduced from the Koch-Tataru result, the proof we present here (Section \ref{s:gexist}) is short, `elementary' and essentially relies only on the decay of heat flows (Section \ref{s:heat-decay}), and a stochastic representation of the Navier-Stokes equations using particle trajectories (Section \ref{s:sl-intro}, see also \cites{detsns,thesis}).

\section{Notational conventions and description of results}\label{s:notation}
In this section we describe the notational convention we use, and state the main result we prove. Let $L > 0$ be a fixed length scale, and $\I = [0,L]$. We define the H\"older norms and semi-norms on $\I^d$ by
\begin{gather*}
\hseminorm{u}{\alpha} = \sup_{x,y \in \I^d} L^\alpha \frac{|u(x) - u(y)|}{|x-y|^\alpha}\\
\cnorm{u}{k} = \sum_{|m|\leqs k} L^{|m|} \sup_{\I^d} |D^m u|\\
\hnorm{u}{k}{\alpha} = \cnorm{u}{k} + \sum_{|m| = k} L^k \hseminorm{D^m u}{\alpha}
\end{gather*}
where $D^m$ denotes the derivative with respect to the multi index $m$. We let \cspace{k} denote the space of all $k$-times continuously differentiable spatially periodic functions on $\I$, and \holderspace{k}{\alpha} denote the space of all spatially periodic $k+\alpha$ H\"older continuous functions. The spaces \cspace{k} and \holderspace{k}{\alpha} are endowed with the norms $\cnorm{\cdot}{k}$ and $\hnorm{\cdot}{k}{\alpha}$ respectively.

We use $\Ifn$ to denote the identity function on $\R^d$ (or on $\I^d$ depending on the context), and use $\Imatrix$ to denote the identity matrix. The main theorem we prove in this paper is
\begin{theorem}\label{t:norm-decay}
Let $k \geq 1$, $\alpha \in (0,1)$ and $u_0 \in \holderspace{k+1}{\alpha}(\I^d)$ be spatially periodic, divergence free and have mean $0$. Let $R = \frac{L}{\nu}\hnorm{u_0}{k+1}{\alpha}$ be the Reynolds number of the flow. Then $\exists T=T(k, \alpha, d, \frac{1}{L}\hnorm{u_0}{k+1}{\alpha})$ and $R_0 = R_0(k, \alpha, d)$ such that for all $R < R_0$ the solution $u$ of \navierstokes with viscosity $\nu = \frac{L}{R} \hnorm{u_0}{k+1}{\alpha}$, initial data $u_0$ and periodic boundary conditions is in \holderspace{k+1}{\alpha} for time $T$, and satisfies
\begin{equation}\label{e:norm-decay}
\hnorm{u_T}{k+1}{\alpha} \leq \hnorm{u_0}{k+1}{\alpha}
\end{equation}
\end{theorem}

We prove Theorem \ref{t:norm-decay} in Section \ref{s:gexist}. A few remarks are in order.

\begin{remark}\label{r:small-data}
Local existence (Theorem \ref{t:snsexist}) combined with the Theorem \ref{t:norm-decay} immediately show that for given initial data, we can choose $\nu$ large enough so that \navierstokes have time global \holderspace{k+1}{\alpha} solutions. Alternately for fixed viscosity, if $\hnorm{u_0}{k+1}{\alpha}$ is small enough, Theorems \ref{t:norm-decay} and \ref{t:snsexist} again give time global \holderspace{k+1}{\alpha} existence of \navierstokes.
\end{remark}

\begin{remark}
The assumption that $u_0$ has mean $0$ is not restrictive. First note that our boundary conditions imply that $\int u_t$ is conserved in time. Set $\bar{u} = \frac{1}{L^d} \int u_0$ to be the mean velocity. Now if we change to coordinates moving with the mean velocity by letting $u'(x, t) = u(x + \bar{u} t, t) - \bar{u}$, then $u'$ solves \navierstokes with mean $0$ initial data $u_0 - \bar{u}$. Thus the smallness assumption in Remark \ref{r:small-data} is really smallness assumptions on the deviation from the mean velocity.
\end{remark}

\begin{remark}
Theorem \ref{t:norm-decay} shows that for some time $T$, equation \eqref{e:norm-decay} holds. Unfortunately our proof does not show that $\hnorm{u_t}{k+1}{\alpha}$ is decreasing in time.
\end{remark}

\section{The stochastic Lagrangian formulation}\label{s:sl-intro}
The Kolmogorov forward equation (or Feynman-Kac formula) \cites{friedman,karatzas} have been extensively used to represent solutions of linear parabolic PDE's as the average of a stochastic process. In this section we briefly describe here a different approach developed in \cites{detsns,thesis,sentropy}, which we use to provide a representation of the Navier-Stokes equations based on noisy particle paths.

Let $u:\R^d \times [0,\infty) \to \R^d$ be some given (time dependent) vector field, and $\theta$ a solution to the heat equation
\begin{equation}\label{e:heat}
\del_t \theta + (u \cdot \grad) \theta - \nu \lap \theta = 0
\end{equation}
with initial data $\theta_0$. We impose either periodic or decay at infinity boundary conditions on $\theta$.

We express $\theta$ as the expected value of a stochastic process as follows: Let $W$ be a $d$ dimensional Wiener process, and let $X:\R^d \to \R^d$ be a solution to the SDE
$$
dX = u \,dt + \sqrt{2 \nu} \,dW
$$
with initial data $X_0(a) = a$. Standard theory\footnote{See also \cites{detsns,thesis,sentropy} for an elementary proof for flows of the type we consider here}~\cite{kunita} shows that the flow $X$ is a homeomorphism, and as spatially differentiable as $u$. We let $A_t$ denote the spatial inverse of the flow map $X_t$.

\begin{proposition}\label{p:rcharacteristics}
If $u \in \cspace{1}$, $\theta_0 \in \cspace{2}$ then the unique solution $\theta$ of \eqref{e:heat} with initial data $\theta_0$ and either periodic or decay at infinity boundary conditions is given by
\begin{equation}\label{e:sl-rep}
\theta_t = \E \theta_0( A_t )
\end{equation}
where $\E$ denotes the expected value with respect to the Wiener measure.
\end{proposition}

Note that if $\nu = 0$, then Proposition \ref{p:rcharacteristics} is nothing but the method of characteristics. If $\nu > 0$, this can be interpreted as solving along random characteristics, and then averaging. Notice also that the Wiener process $\sqrt{2\nu} W_t$ is the natural one to consider here, as it's generator is $\nu\lap$.

The reason we use the representation \eqref{e:sl-rep} and not the Kolmogorov forward equation is because the Kolmogorov forward equation in it's natural setting involves final conditions, and not initial conditions. Thus the standard method employed by probabilists is to make a $t = T - s$ substitution~\cite{freidlin}. The process obtained in this manner will have the same one dimensional distribution as the process $A_t$ above, however spatial covariances and gradients of the two processes will in general be different. The stochastic representation of the Navier-Stokes equations we describe below involves spatial gradients of the flow map $A$, and for this reason our representation will not be valid if we use the Feynman Kac formula.

We now use Proposition \ref{p:rcharacteristics} to represent the solution to the Navier-Stokes equations as the expected value of a system that is nonlinear in the sense of McKean. The essential idea is to find a representation of the Euler equations involving particle trajectories \cites{ele}, and then add noise and average as in Proposition \ref{p:rcharacteristics} (as opposed to attempting to use the Kolmogorov forward equation).

\begin{theorem}\label{t:sl-ns}
Let $\nu > 0$, $W$ be an $n$-dimensional Wiener process, $k\geq1$ and $u_0 \in \holderspace{k+1}{\alpha}$ be a given deterministic divergence free vector field. Let the pair $u$, $X$ satisfy the stochastic system
\begin{align}
\label{e:X-def} dX_t &= u_t \,dt + \sqrt{2 \nu} \,dW_t\\
\label{e:A-def} A_t &= X_t\inv\\
\label{e:u-def} u_t &= \E \lhp\left[ (\gradt A_t) \, (u_0 \circ A_t) \right]
\intertext{with initial data}
\label{e:X-idata} X(a,0) &= a.
\end{align}
We impose boundary conditions by requiring $u$ and $X-I$ are either spatially periodic, or decay at infinity. Then $u$ satisfies the incompressible Navier-Stokes equations \navierstokes with initial data $u_0$.
\end{theorem}

Here $\lhp$ in equation \eqref{e:u-def} denotes the Leray-Hodge projection onto divergence free vector fields \cites{chorin}. We remark that \eqref{e:u-def} is algebraically equivalent to
\begin{align}
\label{e:vorticity} \omega_t &= \E [(\grad X_t) u_0] \circ A_t\\
\label{e:biot-savart} u_t &= - \lap\inv \curl \omega
\end{align}
and \eqref{e:u-def} can be replaced with \eqref{e:vorticity}--\eqref{e:biot-savart} in Theorem \ref{t:sl-ns}. Note that \eqref{e:biot-savart} is exactly the Biot-Savart law. When $\nu = 0$, equation \eqref{e:vorticity} reduces to the well known vorticity transport for the Euler equations \cites{chorin}, and in this case \snseqns (or equivalently the system \eqref{e:X-def}, \eqref{e:A-def}, \eqref{e:X-idata}--\eqref{e:biot-savart}) are exactly a Lagrangian formulation of the Euler equations \cite{ele}.

We do not prove Proposition \ref{p:rcharacteristics} or Theorem \ref{t:sl-ns} here, and we refer the reader to \cites{detsns,thesis} instead. For a generalization of Proposition \ref{p:rcharacteristics} where the diffusion matrix is not spatially constant we refer the reader to \cites{thesis,sentropy}.

\section{Decay of heat flows}\label{s:heat-decay}
In this section we prove a decay estimate for solutions to the heat equation with an incompressible drift. Our first estimate is an $L^\infty \to L^\infty$ estimate that is independent of the drift. A more general $L^1 \to L^\infty$ version of this estimate appeared for example in \cite{mixing} and \cite{quenching}. We provide a proof that follows the proof in \cite{mixing} and keeps track of the dependence of the constants on viscosity and our length scale $L$.
\begin{lemma}\label{l:heat-decay}
Let $u \in C^1( [0,T], \I^d )$ be divergence free, and $\theta$ be a solution to the equation \eqref{e:heat} with initial data $\theta_0$. If $\theta_0$ is spatially periodic, mean $0$, and the dimension $d \geq 3$, then there exists an constant $c = c(d)$ such that
\begin{equation*}
\norm{\theta_t}_\infty \leq \frac{c L^d}{ (\nu t)^{d/2}} \norm{\theta_0}_\infty
\end{equation*}
\end{lemma}
\begin{proof}
Let $\varphi$ be mean zero and periodic, $p \geq \frac{d+2}{4}$ and $c = c(d, p)$ be a constant that changes from line to line. Then the H\"older, Poincar\'e and Sobolev in inequalities give
\begin{align*}
\int \varphi^2 &= \int \varphi^{1/p} \varphi^{(2p-1)/p}\\
    &\leq \lpnorm{\varphi}{1}^{1/p} \lpnorm{\varphi}{(2p-1)/(p-1)}^{(2p-1)/p}\\
    &\leq c L^{(4p - d - 2)/2p} \lpnorm{\varphi}{1}^{1/p} \lpnorm{\grad\varphi}{2}^{(2p-1)/p}
\end{align*}
If we set $q = \frac{2}{2p - 1}$ this gives
\begin{equation*}
\lpnorm{\grad\varphi}{2}^2 \geq c L^{(qd-4)/2} \lpnorm{\varphi}{2}^{2+q} \lpnorm{\varphi}{1}^{-q}
\end{equation*}

Now let $\theta'$ and $\theta''$ to be solutions of \eqref{e:heat} with initial conditions $\theta_0^-$ and $\theta_0^+$ respectively. Integrating \eqref{e:heat} immediately shows that $\int \theta'$ and $\int \theta''$ are conserved. Since $\theta'$ and $\theta''$ are of constant sign, this means that $\lpnorm{\theta'}{1}$ and $\lpnorm{\theta''}{1}$ are conserved in time. Finally, the maximum principle implies that $\theta' \leq \theta \leq \theta''$, and hence $\lpnorm{\theta}{1}$ is nondecreasing in time.

Thus multiplying \eqref{e:heat} by $\theta$ and integrating over $\I^d$ gives
$$
\del_t \lpnorm{\theta}{2}^2 = - 2 \nu \lpnorm{\grad \theta}{2}^2 \leq -c \nu L^{(qd - 4)/2} \lpnorm{\theta}{1}^{-q} \lpnorm{\theta}{2}^{2+q} \leq -c \nu L^{(qd - 4)/2} \lpnorm{\theta_0}{1}^{-q} \lpnorm{\theta}{2}^{2+q}.
$$
Dividing by $\lpnorm{\theta}{2}^{2+q}$ and integrating in time gives
$$
\lpnorm{\theta}{2} \leq c \frac{L^{2/q - d/2}}{(\nu t)^{1/q}} \lpnorm{\theta_0}{1}
$$
Let $\mathcal{P}_t(u)$ be the solution operator of \eqref{e:heat}. The above estimate shows
$$
\norm{\mathcal{P}_t(u)}_{L^1 \to L^2} \leq c \frac{L^{2/q - d/2}}{(\nu t)^{1/q}}
$$
Since $u$ is divergence free the dual operator $\mathcal{P}^*(u) = \mathcal{P}(-u)$, and hence satisfies the same bound. Thus
\begin{align*}
\norm{\mathcal{P}_{2t}}_{L^1\to L^\infty} &\leq \norm{\mathcal{P}_t}_{L^1\to L^2} \norm{\mathcal{P}_t}_{L^2 \to L^\infty}\\
    &= \norm{\mathcal{P}_t}_{L^1\to L^2} \norm{\mathcal{P}_t^*}_{(L^\infty)^* \to L^2}\\
    &\leq \norm{\mathcal{P}_t}_{L^1\to L^2} \norm{\mathcal{P}_t^*}_{L^1 \to L^2}\\
    &\leq c \frac{L^{4/q - d}}{(\nu t)^{2/q}}.
\end{align*}
Hence
\begin{equation*}
\lpnorm{\theta}{\infty} \leq c \frac{L^{4/q - d}}{(\nu t)^{2/q}} \lpnorm{\theta_0}{1} \leq c \frac{L^{4/q}}{(\nu t)^{2/q}} \lpnorm{\theta_0}{\infty}
\end{equation*}
Finally, $p \geq \frac{d+2}{4}$ is the same as $q \leq \frac{4}{d}$, and choosing $q = \frac{4}{d}$ concludes the proof.
\end{proof}
\begin{remark}
When $d = 2$, $p \geq \frac{d+2}{4}$ needs to be replaced with $p > \frac{d+2}{4}$, and hence the above proof will show that for any $\epsilon > 0$,
\begin{equation*}
\norm{\theta_t}_\infty \leq \frac{c_\epsilon L^{d+\epsilon}}{ (\nu t)^{(d+\epsilon)/2}} \norm{\theta_0}_\infty
\end{equation*}
\end{remark}

The formulation of the Navier-Stokes equations \snseqns involves recovering the velocity from the inverse flow map via a singular integral operator (either Biot-Savart or Leray-Hodge). The unboundedness of singular integrals on $L^\infty$ (see \cites{stein}) causes Lemma \ref{l:heat-decay} to be insufficient for our purposes. We now extend Lemma \ref{l:heat-decay} to H\"older spaces for use in our global existence proof. Using the stochastic flows from \cites{detsns,thesis} we obtain a H\"older estimate for solutions of \eqref{e:heat} in an elementary manner.

The usual PDE methods \cite{krylov} provide H\"older estimates that grow exponentially in time. The estimate we provide here will in general grow exponentially in time, however decays in time when the viscosity is large, or drift $U$ is small.

\begin{lemma}\label{l:heat-holder}
Let $d \geq 3$ and $u \in \holderspace{k+1}{\alpha}( [0,T], \I^d )$ be divergence free and define $U$ by
\begin{equation}\label{e:ubound-def}
U = \sup_{t \in [0,T]} \hnorm{u_t}{k+1}{\alpha}.
\end{equation}
Let $\theta_0 \in C(\I^d)$ have mean $0$, and $\theta$ satisfy equation \eqref{e:heat} with initial data $\theta_0$. Then there exists $T' = T'( \frac{U}{L}, d, k, \alpha )$ and a constant $c = c( \frac{UT}{L}, d, k, \alpha )$ such that
\begin{equation*}
\hnorm{\theta}{k+1}{\alpha} \leq c \left(\frac{L^d}{(\nu t)^{d/2}} + \left(\frac{U t}{L}\right)^\alpha \right) \hnorm{\theta_0}{k+1}{\alpha}
\end{equation*}
holds for all $t \in [0, T']$. If $d = 2$, the above estimate is still true if we replace $\frac{L^{d}}{(\nu t)^{d/2}}$ with $\frac{c_\epsilon L^{d+\epsilon}}{(\nu t)^{(d+\epsilon)/2}}$ for any $\epsilon > 0$.
\end{lemma}
\begin{remark}
Note that the growth term is independent of the viscosity, and the decay term is independent of the drift $u$.
\end{remark}
\begin{proof}
We present the proof for $d \geq 3$. The $d=2$ case will then follow by replacing $d$ with $d+\epsilon$. Define $X, A$ by equations \eqref{e:X-def} and \eqref{e:A-def} respectively. From \cites{detsns,thesis} and uniqueness of strong solutions to \eqref{e:heat} we know
$$
\theta_t = \E \theta_0 \circ A_t.
$$
Let $\ell = A - \Ifn$ be the Lagrangian displacement. First notice that if $f \in \holderspace{0}{\alpha}$ then Lemma \ref{l:grad-lambda} shows
\begin{equation}
\label{e:hseminorm-foa} \hseminorm{f \circ A_t}{\alpha} \leq c \hseminorm{f}{\alpha} \left(\frac{U t}{L}\right)^\alpha	\qquad\as
\end{equation}
Now, let $m$ a multi index with $1 \leq \abs{m} \leq k$. We note that $D^m \theta_t$ is a sum of terms of the form
\begin{equation}\label{e:Dm-terms}
D^n \theta_0\at{A_t} \prod_{1\leq i \leq\abs{n}} D^{n_i} \ell_t \qquad\text{and}\qquad D^n \theta_0\at{A_t}
\end{equation}
where $n_i$'s are multi indices with $\abs{n_i}\geq 1$ and $\abs{n} + \sum_i \abs{n_i} = \abs{m}$. By Proposition \ref{p:rcharacteristics} we know that $\E D^n \theta_0 \at{A_t}$ satisfies \eqref{e:heat} with initial data $D^n \theta_0$, and hence by Lemma \ref{l:heat-decay} we know
\begin{equation*}
\norm{\E [D^n \theta_0] \circ A_t}_{L^\infty} \leq \frac{c L^d}{(\nu t)^{d/2}}.
\end{equation*}
Thus using Lemma \ref{l:grad-lambda}, inequality \eqref{e:hseminorm-foa} we have
\begin{equation}\label{e:hnorm-Dntheta}
\hnorm{\E [D^n \theta_0] \circ A_t}{0}{\alpha} \leq c \left( \frac{L^d}{(\nu t)^{d/2}} + \left(\frac{U t}{L}\right)^\alpha \right).
\end{equation}
Using \eqref{e:hnorm-Dntheta} and Lemma \ref{l:grad-lambda}, we bound the remaining terms of \eqref{e:Dm-terms}, concluding the proof.
\end{proof}
\section{Global existence}\label{s:gexist}
In this section we prove Theorem \ref{t:norm-decay}. We start with a Lemma involving bounds for the Leray-Hodge projection.
\begin{lemma}\label{l:lhp-bounds}
Let $k \geq 1$, and $A, v \in \holderspace{k+1}{\alpha}$ be such that $\grad A$, $v$ are spatially periodic. There exists a constant $c = c(d, \alpha)$ such that
$$
\hnorm{\lhp[ (\gradt A) v ]}{k+1}{\alpha} \leq c \hnorm{\grad A}{k}{\alpha} \hnorm{v}{k+1}{\alpha}
$$
\end{lemma}
\begin{proof}
Since $\lhp$ vanishes on gradients, we can `integrate by parts' to avoid the loss of derivatives. Note
\begin{equation*}
\lhp[ (\gradt u) v ] = \lhp[ \grad( u \cdot v) - (\gradt v) u ] = - \lhp[ (\gradt v) u ]
\end{equation*}
for any $u, v \in \cspace{1}$. Thus we have
\begin{equation*}
\del_i \lhp[ (\gradt A) v ] = \lhp[ (\gradt A) \del_i v ] - \lhp[ (\gradt v) \del_i A ]
\end{equation*}
Since $\lhp$ is Calderon-Zygmund singular integral operator, it is bounded on H\"older spaces \cites{stein,cz-periodic}. Finally note that the right hand side only depends on first derivatives of $A$ and $v$, and the lemma follows by taking H\"older norms.
\end{proof}

We now prove Theorem \ref{t:norm-decay}. We restate it here for the readers convenience.
\begin{theorem*}[\ref{t:norm-decay}]
Let $k \geq 1$, $\alpha \in (0,1)$ and $u_0 \in \holderspace{k+1}{\alpha}(\I^d)$ be spatially periodic, divergence free and have mean $0$. Let $R = \frac{L}{\nu}\hnorm{u_0}{k+1}{\alpha}$ be the Reynolds number of the flow. Then $\exists T=T(k, \alpha, d, \frac{1}{L}\hnorm{u_0}{k+1}{\alpha})$ and $R_0 = R_0(k, \alpha, d)$ such that for all $R < R_0$ the solution $u$ of \navierstokes with viscosity $\nu = \frac{L}{R} \hnorm{u_0}{k+1}{\alpha}$, initial data $u_0$ and periodic boundary conditions is in \holderspace{k+1}{\alpha} for time $T$, and satisfies
\begin{equation}\tag{\ref{e:norm-decay}}
\hnorm{u_T}{k+1}{\alpha} \leq \hnorm{u_0}{k+1}{\alpha}
\end{equation}
\end{theorem*}
\begin{proof}
We assume that $d \geq 3$. The $d = 2$ case follows similarly by replacing $d$ with $d+\epsilon$. Let $C, \delta$ be the dimensional constants in Theorem \ref{t:snsexist}. We let $U = C \hnorm{u_0}{k+1}{\alpha}$, and choose $T$ such that $\frac{UT}{L} < \delta$. By Theorem \ref{t:snsexist}, there exist a pair of \holderspace{k+1}{\alpha} functions $X, u: [0,T] \to \I^d$ which are the unique (strong) solution to \snseqns. Recall that $\hnorm{u_t}{k+1}{\alpha} \leq U$ for all $t \in [0,T]$.

From equation \eqref{e:u-def} we see
$$
u_t = \E \lhp u_0 \circ A_t + \E \lhp \left( \gradt \ell \right) u_0 \circ A_t.
$$
Let $c = c(k, \alpha, d)$ be a constant that changes from line to line. Applying Lemma \ref{l:lhp-bounds} and Lemma \ref{l:grad-lambda} to the second term we have
$$
\hnorm{u_t}{k+1}{\alpha} \leq c\hnorm{ \E u_0 \circ A_t }{k+1}{\alpha} + c\E \hnorm{\gradt \ell_t}{k}{\alpha} \hnorm{u_0 \circ A_t}{k+1}{\alpha}
$$
and hence by Lemma \ref{l:heat-holder} we have
$$
\hnorm{u_t}{k+1}{\alpha} \leq c \left(\frac{L^{d}}{(\nu t)^{d/2}} + \left(\frac{U t}{L}\right)^\alpha \right) \hnorm{u_0}{k+1}{\alpha}.
$$
Minimizing $\frac{L^{d}}{(\nu t)^{d/2}} + \bigl(\frac{U t}{L}\bigr)^\alpha$ in time shows that the minimum value is attained at $t_0 = \frac{c L}{U} R^{d/(2\alpha + d)}$, and the minimum value is $c R^{\alpha d/ (2\alpha + d)}$. Thus we can choose $R$ small enough to ensure $t_0 < T$ and equation \eqref{e:norm-decay} is satisfied.
\end{proof}
\appendix
\section{Bounds for the Lagrangian displacement}\label{s:sl-bounds}
In this section, we prove bounds on $\hnorm{ \grad X - \Imatrix}{k}{\alpha}$. The estimates proved here are elementary, and are taken directly from \cites{sperturb,thesis}. We reproduce them here for completeness and the readers convenience.

We remark that the estimates provided here were used in \cites{sperturb,thesis} to prove local existence for the system \snseqns. As the local existence proof is a little lengthier, we do not reproduce it here.

\begin{lemma}\label{l:inverse-bound}
Let $X$ be a Banach algebra. If $x \in X$ is such that $\norm{x} \leqs \rho < 1$ then $1 + x$ is invertible and $\norm{(1 + x)\inv} \leqs \frac{1}{1 - \rho}$. Further if in addition $\norm{y} \leqs \rho$ then
$$\norm{(1+x)\inv - (1+y)\inv} \leqs \frac{1}{(1 - \rho)^2} \norm{x-y}$$
\end{lemma}
\begin{proof}
The first part of the Lemma follows immediately from the identity $(1+x)\inv = \sum (-x)^n$. The second part follows from the first part and the identity
\begin{equation*}
(1+x)\inv - (1+y)\inv = (1+x)\inv (y-x) (1+y)\inv.\qedhere
\end{equation*}
\end{proof}
\begin{lemma}\label{l:holder-ineq}
If $k \geqs 1$, then there exists a constant $c = c(k, \alpha)$ such that
\begin{gather*}
\hnorm{f \circ g}{k}{\alpha} \leqs c \hnorm{f}{k}{\alpha}\left(1 + \hnorm{\grad{g}}{k-1}{\alpha}\right)^{k+\alpha}\\
\shortintertext{and}
\begin{multlined}
\hnorm{ f\circ g_1 - f \circ g_2}{k}{\alpha} \leqs c \hnorm{\grad f}{k}{\alpha} \left( 1 + \hnorm{\grad g_1}{k-1}{\alpha} + \hnorm{\grad g_2}{k-1}{\alpha} \right)^{k+1} \cdot\\
\cdot \hnorm{g_1 - g_2}{k}{\alpha}.
\end{multlined}
\end{gather*}
\end{lemma}

The proof of Lemma \ref{l:holder-ineq} is elementary and not presented here.

\begin{lemma}\label{l:inverse-close}
Let  $X_1, X_2 \in \holderspace{k+1}{\alpha}$ be such that
$$\hnorm{\grad X_1 - \Imatrix}{k}{\alpha} \leqs d < 1 \quad\text{and}\quad \hnorm{\grad X_2 - \Imatrix}{k}{\alpha} \leqs d < 1.$$
Let $A_1$ and $A_2$ be the inverse of $X_1$ and $X_2$ respectively. Then there exists a constant $c = c(k, \alpha, d)$ such that
$$\hnorm{A_1 - A_2}{k}{\alpha} \leqs c \hnorm{X_1 - X_2}{k}{\alpha}$$
\end{lemma}
\begin{proof}
Let $c = c(k, \alpha, d)$ be a constant that changes from line to line (we use this convention implicitly throughout this paper). Note first $\grad A = (\grad X)\inv \circ A$, and hence by Lemma \ref{l:inverse-bound}
$$\cnorm{\grad A}{0} \leqs \cnorm{(\grad X)\inv}{0} \leqs c.$$
Now using Lemma \ref{l:inverse-bound} to bound $\hnorm{ (\grad X)\inv}{0}{\alpha}$ we have
$$\hnorm{\grad A}{0}{\alpha} = \hnorm{ (\grad X)\inv \circ A}{0}{\alpha} \leqs \hnorm{ (\grad X)\inv}{0}{\alpha} \left( 1 + \cnorm{\grad A}{0} \right) \leqs c$$
When $k \geqs 1$, we again bound $\hnorm{(\grad X)\inv}{k}{\alpha}$ by Lemma \ref{l:inverse-bound}. Taking the \holderspace{k}{\alpha} norm of $(\grad X)\inv \circ A$ we have
$$\hnorm{\grad A}{k}{\alpha} \leqs \hnorm{(\grad X)\inv}{k}{\alpha} \left( 1 + \hnorm{\grad{A}}{k-1}{\alpha} \right)^k.$$
So by induction we can bound $\hnorm{\grad A}{k}{\alpha}$ by a constant $c = c(k, \alpha, d)$. The Lemma now follows by applying Lemma \ref{l:holder-ineq} to the identity
\begin{align*}
A_1 - A_2 &= \left( A_1 \circ X_2 - \Ifn \right) \circ A_2\\
    &= \left(A_1 \circ X_2 - A_1 \circ X_1 \right) \circ A_2. \qedhere
\end{align*}
\end{proof}
\begin{lemma}\label{l:grad-lambda}
Let $u \in C([0,T], \holderspace{k+1}{\alpha})$ and $X$ satisfy the SDE \eqref{e:X-def} with initial data \eqref{e:X-idata}. Let $\lambda = X - \Ifn$ and $U = \sup_t \hnorm{u(t)}{k+1}{\alpha}$. Then there exists $T = T(k, \alpha, \frac{U}{L})$ and $c = c(k, \alpha, \frac{UT}{L})$ such that for $t \leq T$
\begin{equation*}
\hnorm{\grad \lambda(t)}{k}{\alpha} \leqs \frac{c U t}{L} \qquad\text{and}\qquad\hnorm{\grad \ell(t)}{k}{\alpha} \leqs \frac{c U t}{L}.
\end{equation*}
hold almost surely.
\end{lemma}
\begin{proof}
From equation \eqref{e:X-def} we have
\begin{alignat}{2}
\nonumber	&& X(x,t) &= x + \int_0^t u( X(x,s), s) \,ds + \sqrt{2\nu} W_t\\
\label{e:gradx-def}	&\implies\quad&	\grad X(t) &= I + \int_0^t (\grad u)\circ X \cdot \grad X.
\end{alignat}
Taking the \cspace{0} norm of equation \eqref{e:gradx-def} and using Gronwall's Lemma we have
\begin{equation*}
\cnorm{\grad \lambda(t)}{0} = \cnorm{ \grad X(t) - I}{0} \leqs e^{U t/L} - 1.
\end{equation*}
Now taking the \holderspace{k}{\alpha} norm in equation \eqref{e:gradx-def} we have
$$ \hnorm{ \grad \lambda(t) }{k}{\alpha} \leqs  c \int_0^t \hnorm{\grad u}{k}{\alpha} \left(1 + \hnorm{\grad \lambda}{k-1}{\alpha} \right)^{k+\alpha} \left( 1 + \hnorm{ \grad \lambda}{k}{\alpha}\right).$$
The bound for $\hnorm{\grad \lambda}{k}{\alpha}$ now follows from the previous two inequalities, induction and Gronwall's Lemma. The bound for $\hnorm{\grad \ell}{k}{\alpha}$ then follows from Lemma \ref{l:inverse-close}.

We draw attention to the fact that the above argument can only bound $\grad \lambda$, and not $\lambda$. Fortunately, our results only rely on a bound of $\grad \lambda$.
\end{proof}

We conclude this appendix by stating a slightly modified version theorem which appeared in \cite{sperturb}. The only modification we make is that we trace the dependence of the constants in \cite{sperturb} to dimension less quantities, instead of absolute ones. The proof that appeared in \cite{sperturb} goes through verbatim.

\begin{theorem}\label{t:snsexist}
Let $k \geqs 1$ and $u_0 \in \holderspace{k+1}{\alpha}$ be divergence free. There exists absolute constants $\delta = \delta( k, \alpha, d )$ and $C = C( k, \alpha, d)$ such that for $U = C \hnorm{u_0}{k+1}{\alpha}$, and any $T$ such that $\frac{U T}{L} < \delta$ there exist a pair of functions a pair of functions $\lambda, u \in C([0,T], \holderspace{k+1}{\alpha})$ such that $u$ and $X = \Ifn + \lambda$ satisfy the system \snseqns. Further for all $t \in [0,T]$ we have $\hnorm{u_t}{k+1}{\alpha} \leq U$.
\end{theorem}

\section*{Acknowledgement}

I would like to thank Lenya Ryzhik for pointing out \cite{mixing} which provides the elegant proof of Lemma \ref{l:heat-decay} reproduced here.

\end{document}